\documentclass[12pt]{article}
\usepackage{amsmath}

\newcommand{\alinea}{\hspace{1em}}

\newtheorem{thm}{\bf \alinea Theorem}

\newtheorem{defi}{\bf \alinea Definition}

\newcommand{\Pf}{\nd \alinea Proof. \alinea \em}
\newcommand{\qd}{\hfill qed}
\newcommand{\nd}{\noindent}

\newcommand{\z}{\zeta}

\newcommand{\om}{\omega}

\begin{document}

\title{Twisted Klein curves modulo 2}

\author{Iwan M. Duursma}

\date{version 30Jan2001}

\maketitle

\begin{abstract}
We give an explicit description of all 168 quartic curves over the field 
of two elements that are isomorphic to the Klein curve over an algebraic 
extension. 
Some of the curves have been known for their small class number \cite{MQ},
\cite{LMQ}, others for attaining the maximal number of rational points \cite{S}, \cite{PSW}.
\end{abstract}

The Klein curve is the famous curve
\[
K : X^3Y+Y^3Z+Z^3X = 0
\]
The equation goes back to F. Klein who derived it as a model for the 
modular curve $X(7)$ \cite{K}. In characteristic zero, the curve is the unique curve of
genus three, up to isomorphism, with the maximal number of 168 automorphisms.
Various aspects of the curve are described in the recent papers \cite{JOS}, \cite{PD}, \cite{AS},
in \cite{E}, and in other contributions to \cite{EW}. 
For further and older references dealing with the many fascinating combinatorial and 
arithmetic properties of the curve, we refer to these publications. \\

In this paper, we investigate the different models of the Klein curve that 
exist in characteristic two. First, we present six different inequivalent
models that we found at different places in the literature, followed 
by a brief discussion of different models in characteristic zero.
Then we recall the special behaviour of the bitangents in even characteristic 
and we present the main theorem that classifies all 168 twisted
curves in terms of their generic tangent. 
Finally, we show that this description 
immediately reveals the zeta functions for each of the 168 twisted curves. \\

The results of the paper were obtained while visiting the University of
Puerto Rico at Rio Piedras, September 1995 - May 1996. A first preprint was
distributed at the conference ``Geometry, Algorithms and Arithmetic in the 
Theory of Error Correcting Codes,'' Guadeloupe, April 1-5, 1996. The results
were presented under the title ``Twisted Klein curves'' at the AMS-Benelux 
meeting in Antwerpen, May 22-24, 1996. This version is unchanged from the
January 30, 2001 version. We added two references of related interest: 
\cite{BG} (Example 7.1), kindly pointed out to us by David Joyner, 
and \cite{T}.

\section{Six inequivalent models modulo 2}

The behaviour of the Klein curve is special in all the characteristics that divide 168.
In this paper, we are particularly interested in the different models that 
exist in characteristic two. Theorem \ref{thmP} shows that each model is characterised 
by a transformation $P \in PSL(3,2)$, such that the generic tangent of the curve
at the point $(x,y,z)$ passes through $(x^2,y^2,z^2)P^t$ and $(x^8,y^8,z^8)(P^t)^3$.
Two elements of $PSL(3,2)$ define equivalent curves if and only if they belong
to the same conjugacy class. The conjugacy classes of $PSL(3,2)$ are of size 
1 (identity), 24 (elements of order 7 with even trace), 24 (order 7, odd trace), 
56 (order 3), 42 (order 4), and 21 (order 2).
Up to equivalence we find six different models. For each of the six conjugacy classes, 
the corresponding
curve has $0$ (identity), $7$ (order 7, even trace), $0$ (order 7, odd trace), $3$ (order 3), 
$2$ (order 4), or $4$ (order 2) rational points over the binary field. \\

The Klein curve modulo 2 represents a curve with $3$ rational points.
In \cite{K}, Klein gives a different model
\begin{align*} 
&K^\prime : \frac{1}{21 \root{3}\of{7}} (X^4+21X^2YZ-147Y^2Z^2+49XY^3+49XZ^3) = 0
\end{align*}
which modulo 2 represents a curve with $4$ rational points (all points not on the line $Y=Z$).
L.E. Dickson, in his paper ``Quartic curves modulo 2'' \cite{DQ},
classifies up to equivalence {\em all} the quartics that have $0,7,6,$ or $5$ rational
points (for the last case, he finds twenty five equivalence classes and 
``Quartics with 1, 2, 3, or 4 real points have not been treated since they would 
probably not present sufficient novelty to compensate for the increased length
of the investigation''). The paper provides us with three more inequivalent models, 
\begin{align*}
&\alpha : {X}^{4}+{Y}^{4}+{Z}^{4}+{X}^{2}{Y}^{2}+{Y}^{2}{Z}^{2}+{Z}^{2}{X}^{2}+
{X}^{2}YZ+X{Y}^{2}Z+XY{Z}^{2} = 0 \\
&A : {X}^{3}Y+{X}^{2}{Y}^{2}+X{Z}^{3}+{X}^{2}{Z}^{2}+{Y}^{3}Z+Y{Z}^{3} = 0 \\
&\gamma_{1,0} : {X}^{4}+{Y}^{4}+{Z}^{4}+{X}^{3}Y+X{Z}^{3}+{Y}^{3}Z+XY{Z}^{2} = 0
\end{align*}
The model $\alpha$ is the unique model invariant under the full group $PSL(3,2)$.
The model $A$ has the maximum number 
of rational points for a genus 3 curve. An equivalent model appears in \cite{S},
\cite{PSW}. We will denote it $X_{N=7}$.
\begin{align*}
&X_{N=7} : XY(X+Y)(X+Z)+YZ(Z+Y)Y+ZX(X+Z)Z = 0
\end{align*}
The model $\gamma_{1,0}$, with no rational points, has trivial Jacobian.
An equivalent model appears in the list by Madan and Queen \cite{MQ} of the seven 
congruence function fields which have class number one and genus different from zero. 
We will denote it $X_{h=1}$.
\begin{align*}
&X_{h=1} : Y^4+XY^3+(X^2+XZ)Y^2+(X^3+Z^3)Y+(X^4+XZ^3+Z^4)=0 
\end{align*}
The classes for $A$ and $\gamma_{1,0}$ are of same size and differ by a
translation by $\alpha$. It happens that $X_{N=7} = \gamma_{1,0} + \alpha$, and
$X_{N=7}$ and $\gamma_{1,0}$ have the same seven rational automorphisms. 
For the class with $2$ rational points, we give the model
\begin{align*} 
&X_{N=2} : \,{X}^{4}+{Y}^{4}+{Z}^{4}+{X}^{3}Y+X{Z}^{3}+{Y}^{3}Z+XY{Z}^{2}+ \\
&\qquad \qquad +{X}^{2}{Y}^{2}+{X}^{2}{Z}^{2}+{Y}^{2}{Z}^{2}+Y{Z}^{3} = 0 
\end{align*}
Comparison with $\gamma_{1,0}$ shows that the rational points are those not
on $X(Y+Z)=0$. Properties of the different models are summarized in the Appendix.
  

\section{Models in characteristic zero}

The Klein curve arises in characteristic zero as a model for the elliptic
modular curve $X(7)$. In particular this yields that the curve has the
simple group $PSL(2,7)$ of order 168 as automorphism group. Many of the
interesting properties of the Klein curve are directly related to this
group. F. Klein derived the model 
\[
X^3Y+Y^3Z+Z^3X = 0
\]
as follows. 
A plane curve of degree four contains 24 flexpoints, that is points where the
tangent intersects with multiplicity at least three. The group acts faithfully
on the set of flexpoints and each flexpoint has stabilizer of order seven.
There are eight subgroups of order seven, each fixing at most three points.
It follows that each subgroup of order seven fixes three points,
all of them flexpoints.
The subgroup fixing a flexpoint also fixes the fourth point of
intersection on its tangent, which must therefore be a flexpoint with the same
stabilizer.
Three flexpoints connected by their tangents form an oriented triangle, that
Klein called a {\em Wendedreiecke}. The Klein quartic is now
the curve with {\em triangle of inflection tangents} $XYZ=0$. \\

It is convenient to view the action of the automorphism group on the curve
through its action on the eight {\em triangles of inflection tangents}. 
Each triangle has stabilizer
of order 21. Let $\z$ be a primitive seventh root of unity. The stabilizer
of the triangle $XYZ=0$ is generated by
\[
\tau =
\left( \begin{array}{ccc}
\z &0    &0 \\
0  &\z^4 &0 \\
0  &0    &\z^2   \\
\end{array} \right)
\]
and the cyclic permutation $\rho$ of the coordinates $X, Y, Z$. It has 
14 elements of order three. The full group has 56 elements of order three
and each automorphism of order three appears in the stabilizer
of two triangles. Say the other triangle stabilized by the cyclic permutation
$\rho$ has a vertex $(a:b:c)$. Klein showed that
\[
\sigma =
\left( \begin{array}{ccc}
a &b &c \\
b &c &a \\
c &a &b \\
\end{array} \right)
\]
acts on the curve, for 
\[
a = (\z^5-\z^2)/\sqrt{-7}, \, b = (\z^3-\z^4)/\sqrt{-7}, \, c = (\z^6-\z)/\sqrt{-7},
\]
with $\sqrt{-7} = \z+\z^2-\z^3+\z^4-\z^5-\z^6$.
The automorphism $\sigma$ is of order two and interchanges the two triangles.
Indeed, $\sigma$ is of order two only if $f=(x-a)(x-b)(x-c)$ is of the form 
$x^3 \pm x^2 - r$, and it acts on the Klein curve only if $r=0,\pm 1/7$. \\

The full automorphism group $PSL(2,7)$ of the Klein curve is visible from its
action on the eight {\em triangles of inflection tangents}. 
Let $\Delta_\infty$ be the triangle $XYZ=0$, 
let $\Delta_0=\sigma\Delta_\infty$ be the other real triangle, and let 
$\Delta_i = \tau^i \Delta_0$, for $i=0,1,\ldots,6$. 
The action of $\sigma$ and $\tau$ on the triangles is described by the 
permutations $(1\,6)(2\,3)(4\,5)(0\,\infty)$ and $(1\,2\,3\,4\,5\,6\,0)$ respectively.
And an isomorphism to $PSL(2,7)$ is given by
\[
\tau \longmapsto \left( \begin{array}{cc}1 &1 \\ 0 &1 \\ \end{array} \right) \quad
\sigma \longmapsto \left( \begin{array}{cc}0 &-1 \\ 1 &0 \\ \end{array} \right)
\]
Burnside pointed out that $\tau$ and $\sigma$ generate the automorphism group. 
Klein had originally included the cyclic permutation $\rho$ in the generating set. \\

The two {\em triangles of inflection tangents} that are left invariant by the cyclic 
permutation $\rho$ are
the only two of the eight triangles with real vertices. If we change the orientation
in the triangle $XYZ=0$, we obtain the different model
\[
X^3Z+Y^3X+Z^3Y = 0.
\]
If we change the orientation in the other real triangle, we obtain the different model
\begin{multline*}
O_4 : \frac{1}{7}( (X^4+Y^4+Z^4)-3(X^2Y^2+Y^2Z^2+Z^2X^2) \\
+3(X^2YZ+XY^2Z+XYZ^2)+6(X^3Z+Y^3X+Z^3Y)) = 0.
\end{multline*}
This model was derived independently as Equation (1.22) in \cite{E}.
After reduction modulo $p=2$, the model agrees with the invariant model $\alpha$
and has all automorphisms defined over the base field.
Ciani \cite{C} gives a model over $Q(\sqrt{-7})$  
\[
C_4 : X^4+Y^4+Z^4+3\om(Y^2Z^2+Z^2X^2+X^2Y^2) = 0,~~~~\om^2+\om+2=0.
\]
It coincides with Equation (1.11) in \cite{E}.
The following model (\cite{AS}) has the $S_3$ symmetry of $C_4$ and the same binary reduction as $O_4$,
\[
A_4 : 7 s_1 (s_1^3+s_3) - (2s_1^2+s_2)^2 = 0.
\]
All of the last three curves have their automorphisms defined over $Q(\sqrt{-7})$.
The six real flexpoints of the last curve have $s_1 = -1, s_2 = -2, s_3 = +1$, i.e. 
have coordinates $a,b,c$ such that, for $\theta = 2\pi/7$,
\begin{multline*}
(x-a)(x-b)(x-c) = x^3+x^2-2x-1 = \\
= (x-2\cos(\theta))(x-2\cos(2\theta))(x-2\cos(4\theta))
\end{multline*}


\section{Bitangents in even characteristic}

In characteristic zero, a non-singular quartic curve admits 28 bitangents. 
In even characteristic, the number of bitangents is at most $7$. In that case, 
the bitangents define a Fano plane. All non-singular quartic curves with 
seven bitangents are in the same isomorphism class. The unique curve with
bitangents $X, Y, Z, X+Y, Y+Z, X+Y+Z, X+Z$ is the invariant curve
\[
\alpha : X^4+Y^4+Z^4+Y^2Z^2+Z^2X^2+X^2Y^2+X^2YZ+Y^2ZX+Z^2XY = 0
\]
Uniqueness follows for example by considering that the bitangents must touch the 
curve in the fourteen points of $P^2(F_4)$ not in $P^2(F_2)$. 
Other curves with seven bitangents will have their bitangents defined over some
extension of $F_2$. In those cases, the equation of the curve can be recovered 
from the bitangents by substitution in $\alpha$ of three non-concurrent bitangents. \\

The Klein curve has bitangents
\[
K : a^2X + aY + a^4Z = 0, \quad a^8+a=0,~~a \neq 0. 
\]
It happens that the models chosen by L.E. Dickson have a similar
description of the bitangents. The model $A$ has bitangents
\[
A : b^2X + bY + b^4Z = 0, \quad b^8+b^2+b=0,~~b \neq 0.
\]
The model $\gamma_{1,0}$ has bitangents
\[
\gamma_{1,0} : c^2X + cY + c^4Z = 0, \quad c^8+c^4+c=0,~~c \neq 0.
\]
If we choose different but equivalent models, 
the bitangents change but are still parametrized by the same values.
The curve $X_{N=7}$, which is equivalent to $A$, has bitangents 
\[
X_{N=7} : bX + b^{16}Y + b^8Z = 0, \quad b^8+b^2+b=0,~~b \neq 0.
\]
The curve $X_{h=1}$, which is equivalent to $\gamma_{1,0}$, has bitangents
\[
X_{h=1} : cX + c^8Y + c^2Z = 0, \quad c^8+c^4+c=0,~~c \neq 0.
\]
We chose the model $X_{N=2}$ such that it has bitangents
\[
X_{N=2} : d^2X + dY + d^4Z = 0, \quad d^8+d^4+d^2+d=0,~~d \neq 0.
\]
For the invariant curve, the bitangents have coefficients, in matrix form,
\[
\alpha : \left [\begin {array}{ccccccc} 1&0&0&1&0&1&1\\\noalign{\medskip}0&1&0&1
&1&1&0\\\noalign{\medskip}0&0&1&0&1&1&1\end {array}\right ]
\]
So that the coefficients appearing in a given row are the zeros of 
$(x^8+x^4)/x$. After reduction, Klein's model $K^\prime$ has bitangents, 
with $F_4 = \{0,1,\omega,\bar \omega\}$, 
\[
K^\prime : \left [\begin {array}{ccccccc} 1&1&0&0&1&0&1 \\\noalign{\medskip}
0&\omega&\bar \omega&\omega&1&1&\bar \omega \\\noalign{\medskip}
0&\bar \omega&\omega&\bar \omega&1&1&\omega\end {array}\right ]
\]
Substituion of $Z=Z+X$ in $K^\prime$ gives the equivalent reduced model
\begin{align*} 
&X_{N=4} : \,X^2YZ+X^2Y^2+X^2Z^2+Y^2Z^2+XY^3+XZ^3+X^3Y+X^3Z = 0 
\end{align*}
with bitangents,
\[
X_{N=4} : \left [\begin {array}{ccccccc} 1&\omega&\omega&\bar \omega&0&1&\bar \omega\\\noalign{\medskip}
0&\omega&\bar \omega&\omega&1&1&\bar \omega\\\noalign{\medskip}
0&\bar \omega&\omega&\bar \omega&1&1&\omega\end {array}\right ]
\]
The coefficients appearing in a given row are the zeros of 
$(x^8+x^2)/x$. The only non-trivial additive polynomial that can not be
used to define coefficients of a Fano plane is $x^8+x^4+x^2$. \\

In all cases, the set of bitangents is defined over $F_2$.
The Frobenius automorphism acts on the set of bitangents. For a careful
choice of the coefficients of the bitangents (such that the linear expressions
for the bitangents are closed under addition), the action is linear
and can be represented by an element of $PSL(3,2)$.

\begin{defi} \label{defR} To a quartic curve over $F_2$ with seven bitangents,
we associate the unique element $R \in PSL(3,2)$ such that the action
of $R$ on the bitangents agrees with the action of the Frobenius automorphism: 
$(a^2,b^2,c^2)^T = R (a,b,c)^T$, for each bitangent $ax+by+cz=0$.
\end{defi}

In the next section, we give the inverse map, which gives for each $R \in PSL(3,2)$ 
the unique quartic curve such that the Frobenius action on the bitangents corresponds 
with the action of $R$.

\section{Description of the 168 twisted curves}

The group $PSL(3,2)$ is the full automorphism group of the curve $\alpha$.
In \cite{DI}, L.E. Dickson gives the modular invariants of
the general linear group over a finite field. In particular, he shows that
$\alpha$ is the smallest projective invariant of $PSL(3,2)$ acting on $F_2[X,Y,Z]$. 
The other generating invariants are of degree six and seven, and the three invariants 
are algebraicly independent 
\cite{DI}, \cite{B}. 
Consider the additive polynomial in the variable $T$ with coefficients in $F_2[X,Y,Z]$
that cancels at $X,Y,Z$,
\[
\prod_{v \in {\langle X,Y,Z \rangle}} (T + v) \\
~=~ T^{8} + I_4(X,Y,Z) T^{4} + I_6(X,Y,Z) T^{2} + I_7(X,Y,Z) T.
\]
The algebra $F_2[I_4,I_6,I_7]$ is the ring of $PSL(3,2)$
invariants in $F_2[X,Y,Z]$. The coefficients follow by solving
\begin{eqnarray*}
\left (
\begin{array}{ccc}
X &X^2 &X^4 \\
Y &Y^2 &Y^4 \\
Z &Z^2 &Z^4
\end{array}
\right )
\left (
\begin{array}{c}
I_7 \\
I_6 \\
I_4
\end{array}
\right )
=
\left (
\begin{array}{c}
X^8 \\
Y^8 \\
Z^8
\end{array}
\right ) .
\end{eqnarray*}
For $I_4$, we have
\begin{eqnarray*}
I_4 (X,Y,Z) =
\det
\left (
\begin{array}{ccc}
X &X^2 &X^8 \\
Y &Y^2 &Y^8 \\
Z &Z^2 &Z^8
\end{array}
\right )
/
\det
\left (
\begin{array}{ccc}
X &X^2 &X^4 \\
Y &Y^2 &Y^4 \\
Z &Z^2 &Z^4
\end{array}
\right ).
\end{eqnarray*}
\nd Thus, $I_4$ is the unique non-trivial plane curve (i.e. different from a line),
for which the line through a point $(x:y:z)$ and its
Frobenius $F(x:y:z)=(x^2:y^2:z^2)$ also contains $F^3(x:y:z)=(x^8:y^8:z^8)$. \\

\nd Let $C$ be a curve isomorphic to $I_4$. Then $C$ can be represented as a plane
quartic. The canonical divisor of a quartic curve is given by the intersection of 
a line with the curve. Therefore the isomorphism preserves collinear points and
can be represented by a linear transformation.  
Let $(x,y,z) \mapsto A(x,y,z) = (x',y',z')$ be an isomorphism from $C$ to $I_4$.
We have that $(x,y,z)$ is a point on $C$ if and only if $(x',y',z')$ is a point
on $I_4$ if and only if
\[
(x',y',z'), F(x',y',z'), F^3(x',y',z')
\]
are collinear if and only if
\[
A(x,y,z), F(A(x,y,z)),  F^3(A(x,y,z))
\]
are collinear if and only if
\[
(x,y,z), (A^{-1}FA)(x,y,z), (A^{-1}FA)^3(x,y,z)
\]
are collinear. We claim that $(A^{-1}FA)(x,y,z) = (PF)(x,y,z)$,
for some $P \in PSL(3,2)$. Let $A^{(2)}$ be the isomorphism $A$ with coefficients 
squared. We have
\[
(A^{-1}FA)(x,y,z) = (A^{-1}A^{(2)}F)(x,y,z)
\]
and $A^{(2)}A^{-1} \in \text{Aut}(I_4) = PSL(3,2)$, say $A^{(2)}A^{-1} = Q$. The
element $Q$ is determined by the curve $C$ only up to conjugacy, for if we replace
$A$ with $SA$, for $S \in \text{Aut}(I_4) = PSL(3,2)$, 
then $(SA)^{(2)}(SA)^{-1} = SQS^{-1}$. On the other hand,
$A^{-1}QA = A^{-1}A^{(2)}$ does not change when $A$ is replaced with $SA$.  
The latter matrix has coefficients in $F_2$, for 
$Q^2 = A^{(2)}A^{-1}Q=QA^{(2)}A^{-1}$, or
\[
A^{-1}QA = (A^{(2)})^{-1}QA^{(2)} = (A^{-1}QA)^{(2)}
\]
and $A^{-1}A^{(2)} = A^{-1}QA$ belongs to $PSL(3,2)$.

\begin{thm} \label{thmP}
A curve is isomorphic to the invariant curve $I_4$ if and only if
it is of the form
\[
\det ( \left( \begin{array}{c} X \\ Y \\ Z \end{array} \right) 
       P\hspace{-1mm}\left( \begin{array}{c} X^2 \\ Y^2 \\ Z^2 \end{array} \right)
       P^3\hspace{-1mm}\left( \begin{array}{c} X^8 \\ Y^8 \\ Z^8 \end{array} \right) ) / 
\det ( \left( \begin{array}{c} X \\ Y \\ Z \end{array} \right) 
       P\hspace{-1mm}\left( \begin{array}{c} X^2 \\ Y^2 \\ Z^2 \end{array} \right)
       P^2\hspace{-1mm}\left( \begin{array}{c} X^4 \\ Y^4 \\ Z^4 \end{array} \right) ) 
\] 
for some $P \in PSL(3,2)$. The rational automorphisms of such a curve
are given by the elements of $PSL(3,2)$ that commute with $P$. The action
of the Frobenius on the bitangents, as in Definition \ref{defR}, 
is given by the matrix $R = P^t$. The tangent at the point $(x,y,z)$
has intersection divisor 
\[
2\,(x,y,z) + (x^2,y^2,z^2)\,P^t + (x^8,y^8,z^8)\,(P^t)^3
\]
In particular, $P \in PSL(3,2)$ is determined by the generic tangent of the curve. \\

\Pf The curve is uniquely determined by the matrix $P = A^{-1}A^{(2)}$,
for any isomorphism $A$ from $C$ to $I_4$. Now $S \in PSL(3,2)$ 
defines an automorphism of
$C$ if and only if $AS$ is an isomorphism from $C$ to $I_4$ if and only if
$(AS)^{-1} (AS)^{(2)} = S^{-1}PS = P$.
Let the columns of the matrix $B$ contain the coefficients of three non concurrent
bitangents. Then $B^{(2)} = RB$. On the other hand, $B^t$ defines an isomorphism
from $C$ to $I_4$. And $P = (B^t)^{-1}(B^t)^{(2)} = R^t$.
\qd
\end{thm}

The general construction of twisted models for a given object with
non-trivial automorphisms is given by elementary galois cohomology.
The case that we consider here seems to be particularly favourable
for an explicit treatment. The invariant curve $I_4$ and all its
automorphisms are defined over the base field, and the equation of 
the curve $I_4$ is a simple relation among Frobenius images of the 
generic point on the curve. \\
 
As illustration of the theorem, we may write for the Klein curve,
\[
{X}^{3}Y+{Y}^{3}Z+{Z}^{3}X ~=~ {\det\hspace{-1mm}\left( 
\begin {array}{ccc} X&{Y}^{2}&{X}^{8}\\\noalign{\medskip}Y&{Z}^
{2}&{Y}^{8}\\\noalign{\medskip}Z&{X}^{2}&{Z}^{8}\end {array}\right)}
\,/\,
{\det\hspace{-1mm}\left( 
\begin {array}{ccc}
X&{Y}^{2}&{Z}^{4}\\\noalign{\medskip}Y&{Z}^
{2}&{X}^{4}\\\noalign{\medskip}Z&{X}^{2}&{Y}^{4}\end {array}\right)}
.
\]







\section{Zeta functions for the twisted curves}
 
The forms $I_4$ and $I_6$ define a pair of dual curves (up to an inseparable
morphism of degree two), and are birationally equivalent. For both curves,
all points on the tangent at $P$ are images of $P$ under the Frobenius
automorphism $F$. And in both cases, the zeta function is determined by the 
intersection divisor of a generic tangent. For the tangent $L_P$ of $I_4$ at $P$,
\[
L_P \cap I_4 = P + P + F(P) + F^3(P),
\]
and the Frobenius eigenvalues are zeros of
\[
(2+t+t^3) = (1+t)(2-t+t^2).
\]
Similarly, for the tangent $L_P$ of $I_6$ at $P$,
\[
L_P \cap I_6 = P + P + P + P + F^2(P) + F^3(P),
\]
and the Frobenius eigenvalues are zeros of
\[
(4+t^2+t^3) = (2+t)(2-t+t^2).
\]
Thus, both $I_4$ and the normalization of $I_6$ have zeta polynomial
\[
L(t) = (1-t+2t^2)^3.
\]
For the twisted models we see, either by elimination of $P^t$ in the
generic tangent (given in Theorem \ref{thmP}) or by using properties of the zeta function, 
that
the Frobenius eigenvalues are reciprocal zeros of
\[
\prod_{\z^m=1} (1-\z t+2 \z^2 t^2)
\]
where $m$ is the order of the matrix $P$. For the relevant values 
$m=2,3,4,7,$ we find the zeta polynomial immediately from the 
factorization of the above product into irreducibles and the known
number of rational points.
Let $z^+=(1-t+2t^2)$ and $z^-=(1+t+2t^2)$.
\[
\begin{array}{lclcl} 
m=2 : & &1+3t^2+4t^4 = z^+ z^-,     & &L(t) = z^+ (z^-)^2. \\
m=4 : & &1-t^4+16t^8 = z^+ z^- z_2, & &L(t) = (z^+)^2 z^-. \\
m=3 : & &1+5t^3+8t^6 = z^+ z_3,     & &L(t) = z^+ z_3. \\
m=7 : & &1+13t^7+128t^{14} = z^+ {z_7}^+ {z_7}^-. 
& &L(t) = {z_7}^+, {z_7}^-  
\end{array}
\]
The two irreducible factors ${z_7}^+, {z_7}^-$ of degree six
represent the two different conjugacy classes for $P$ of order $7$.
\[
\begin{array}{lcl}
\text{even trace:} & &L(t) = {z_7}^+ = 
1+4\,t+9\,{t}^{2}+15\,{t}^{3}+18\,{t}^{4}+16\,{t}^{5}+8\,{t}^{6} \\
\text{odd trace:}  & &L(t) = {z_7}^- =
1-3\,t+2\,{t}^{2}+{t}^{3}+4\,{t}^{4}-12\,{t}^{5}+8\,{t}^{6} 
\end{array}
\]
The factor $z^+$ corresponds, as is well-known \cite{E}, 
to the elliptic curve with complex multiplication by $-7$.
In the cases $m=2,3,4$, it appears as the quotient of the twist 
by a rational automorphism. In the case $m=7$, the quotient gives the 
projective line. 
We give the elliptic curve as the
quotient, in characteristic zero, of the Klein curve 
by the cyclic permutation of the coordinates.
Let $B$ be the divisor sum of the two points in the intersection of the 
Klein curve and the invariant bitangent $X+Y+Z=0$. 
Recall that the Klein curve
has two {triangles of inflection tangents} that are invariant under the 
cyclic permutation. Let $\Delta_\infty$ and $\Delta_0$, respectively,
be the divisor sums of their vertices, respectively. 
Let $x,y,z \in k[X,Y,Z]$ intersect the Klein curve as follows
\[
(x)=2\Delta_\infty+2\Delta_0, ~
(y)=\Delta_\infty+\Delta_0+3B, ~
(z)=4\Delta_\infty.
\]
We obtain
\begin{align*}
&x=X^2Y+Y^2Z+Z^2X+XYZ=((X^3Y+Y^3Z+Z^3X)+s_s^2)/s_1, \\
&y=(X+Y+Z)(XY+YZ+ZX), \quad z=XYZ.
\end{align*}
The morphism $(-x:y:z)$ maps $K$ onto the elliptic curve
\[
E: y^2+5xy=x^3-x^2+7x.
\]
In fact,
\begin{multline*}
y^2z-5xyz+x^3+x^2z+7xz^2 = \\
= (X^3Y+Y^3Z+Z^3X)(X^3Y^2+Y^3Z^2+Z^3X^2).
\end{multline*}
Replacing $x=((X^3Y+Y^3Z+Z^3X)+s_s^2)/s_1$ by $x=s_2^2/s_1$ in the morphism,
we find, as in \cite{E},
\[
y^2z-5xyz+x^3+x^2z+7xz^2 = (X^3Y+Y^3Z+Z^3X)(X^3Z+Y^3X+Z^3Y)
\]
The triangles $\Delta_\infty$ and $\Delta_0$ are mapped to the point at infinity
and the origin respectively. The only ramification is at $B$.

\appendix

\section{Summary of curves}

The curves appearing in this paper are summarized in the following table.
It gives the number of equivalent curves in a class, the structure of the
subgroup of rational automorphisms, and properties of the matrix $P$.
The matrix $P$ for each of the curves follows (as determined by
Theorem \ref{thmP}).
 
\begin{center}
\begin{tabular}{|l|r|c|c|c|} \hline
Curve                     &$\#$  &Rational   &\multicolumn{2}{|c|}{$P$}  \\
                          &      &automorphisms   &order &trace   \\ \hline
$\alpha=I_4$              &1     &$PSL(3,2)$ &1     &1       \\
$A$, $X_{N=7}$            &24    &$Z/7Z$     &7     &0       \\
$\gamma_{1,0}$, $X_{h=1}$ &24    &$Z/7Z$     &7     &1       \\
$K$                       &56    &$Z/3Z$     &3     &0       \\
$X_{N=2}$                 &42    &$Z/4Z$     &4     &1       \\
$K^\prime$, $X_{N=4}$     &21    &$D_4$      &2     &1       \\
\hline
\end{tabular} \\
\end{center}

\[
\begin{array}{lclclcl}
\alpha   &\left [\begin {array}{ccc} 1&0&0\\\noalign{\medskip}0&1&0\\\noalign{\medskip}0&0&1\end {array}\right ] &
A        &\left [\begin {array}{ccc} 0&1&1\\\noalign{\medskip}0&0&1\\\noalign{\medskip}1&0&0\end {array}\right ] &
X_{N=7}  &\left [\begin {array}{ccc} 1&1&0\\\noalign{\medskip}0&1&1\\\noalign{\medskip}1&0&0\end {array}\right ] \\
K        &\left [\begin {array}{ccc} 0&1&0\\\noalign{\medskip}0&0&1\\\noalign{\medskip}1&0&0\end {array}\right ] &
\gamma_{1,0} &\left [\begin {array}{ccc} 0&1&0\\\noalign{\medskip}0&0&1\\\noalign{\medskip}1&0&1
              \end {array}\right ] &
X_{h=1}  &\left [\begin {array}{ccc} 0&0&1\\\noalign{\medskip}0&1&1\\\noalign{\medskip}1&1&0\end {array}\right ] \\
X_{N=2}  &\left [\begin {array}{ccc} 0&1&1\\\noalign{\medskip}0&0&1\\\noalign{\medskip}1&0&1\end {array}\right ] &
K^\prime &\left [\begin {array}{ccc} 1&0&0\\\noalign{\medskip}1&0&1\\\noalign{\medskip}1&1&0\end {array}\right ] &
X_{N=4}  &\left [\begin {array}{ccc} 1&0&0\\\noalign{\medskip}1&0&1\\\noalign{\medskip}1&1&0\end {array}\right ] 
\end{array}
\]

\noindent
Current address: \\
Department of Mathematics \\                                                  
University of Illinois at Urbana-Champaign \\ 
1409 W. Green Street \\
Urbana, IL 61801 \\
E-mail: duursma@math.uiuc.edu


\begin{thebibliography}{9}

\bibitem{B} Benson, D. J.,
{\it Polynomial invariants of finite groups.} 
London Mathematical Society LNS 190. 
Cambridge University Press, Cambridge, 1993.

\bibitem{BG} Brock, B. W. and Granville, A.,
More points than expected on curves over finite field extensions. 
Finite Fields Appl. 7 (2001), no. 1, 70--91.

\bibitem{C} Ciani, E., 
I varii tipi possibili di quartichi piani p\'{\i}u volte omologico-harmoniche.
Rend. Circ. Mat. Palermo 13 (1899), 347--373.

\bibitem{DI} Dickson, L. E., 
A fundamental system of invariants of the general modular linear group 
with a solution of the form problem. 
Trans. Amer. Math. Soc. 12 (1911), no. 1, 75--98. 

\bibitem{DQ} Dickson, L. E., Quartic curves modulo $2$. 
Trans. Amer. Math. Soc. 16 (1915), no. 2, 111--120. 

\bibitem{AS} Duursma, I. M., Monomial embeddings of the Klein curve. 
Combinatorics (Assisi, 1996). Discrete Math. 208/209 (1999), 235--246.

\bibitem{E} Elkies, N. D., The Klein quartic in number theory.
In: {\it The eightfold way}, 51--101, 
Math. Sci. Res. Inst. Publ., 35, Cambridge Univ. Press, Cambridge, 1999. 

\bibitem{JOS} Jeurissen, R. H., van Os, C. H., and Steenbrink, J. H. M.,
The configuration of bitangents of the Klein curve. 
Discrete Math. 132 (1994), no. 1-3, 83--96.

\bibitem{K} Klein, F., Ueber die Transformation siebenter Ordnung der elliptischen Funktionen. 
Math. Ann. 14 (1879), 428--471. (On the order-seven transformation of elliptic functions. 
Translated from the German by Silvio Levy. Math. Sci. Res. Inst. Publ., 35, {\it The eightfold way}, 
287--331, Cambridge Univ. Press, Cambridge, 1999.)

\bibitem{EW} {\it The eightfold way}. The beauty of Klein's quartic curve. Edited by Silvio Levy. 
Mathematical Sciences Research Institute Publications, 35. 
Cambridge University Press, Cambridge, 1999. 

\bibitem{LMQ} Leitzel, J. R. C., Madan, M. L., and Queen, C. S.,
Algebraic function fields with small class number. 
J. Number Theory 7 (1975), 11--27.

\bibitem{MQ} Madan, M. L. and Queen, C. S.,
Algebraic function fields of class number one. 
Acta Arith. 20 (1972), 423--432.

\bibitem{PSW} Pellikaan, R., Shen, B.Z., and vanWee, G.J.M.,
Which linear codes
are algebraic-geometric, IEEE Trans. Inform. Theory, vol.37, pp.583-602, 1991.

\bibitem{PD} Prapavessi, D. T., On the Jacobian of the Klein curve. 
Proc. Amer. Math. Soc. 122 (1994), no. 4, 971--978. 

\bibitem{S} Serre, J.-P., Sur le nombre de points rationels d'une courbe
alg\'{e}brique sur un corps fini, C.R. Acad. Sci. Paris, vol.296,
pp.397-402, 1983; Oeuvres, III, no.128, pp.658-663.

\bibitem{T} Top, J., Curves of genus 3 over small finite fields.
arXiv:math.NT/0301264.

\end{thebibliography}
\end{document}